\documentclass[12pt]{article}

\author{S.~K.~Lando\thanks{
Russian--French Laboratory for Mathematics, Informatics,
and Mathematical Physics of the
Independent University of Moscow and
Institute for System Research, RAS;
partly supported by the Russian Foundation
for Basic Research (project
01-01-00660), NWO-RFBR (project 047.008.005), and INTAS (project 00259),
{\bf e-mail:} lando@lando.mccme.rssi.ru},
D.~Zvonkine\thanks{b\^at. 425,
Universit\'e Paris-Sud, 91400 Orsay. {\bf e-mail:}
dimitri.zvonkine@math.u-psud.fr  \hspace{5cm}
Both authors are partially supported by the RFBR
grant 02-01-22004.}}

\title{Counting ramified coverings and intersection
theory on spaces of rational functions I \\
{\Large (Cohomology of Hurwitz spaces)}}

\date{\today}

\usepackage{amssymb}
\usepackage{theorem}

\def\C{{\mathbb C}}

\def\Aut{{\rm Aut}}
\def\tLL{{\widetilde{LL}}}

\def\P{{\mathbb P}}
\def\CP{{\mathbb C} {\rm P}}

\def\dim{{\rm dim}}

\def\qed{{\hfill $\diamond$}}

\def\L{{\cal L}}
\def\D{{\cal D}}

\def\cH{{\cal H}}
\def\cL{{\cal L}}
\def\cM{{\cal M}}
\def\cO{{\cal O}}
\def\cT{{\cal T}}
\def\tcT{{\widetilde {\cal T}}}
\def\ocM{{\overline \cM}}
\def\tPsi{{\widetilde \Psi}}
\def\tSigma{{\widetilde \Sigma}}
\def\tPcH{{\widetilde {\P\cH}}}
\def\tPSigma{{\widetilde{\P\Sigma}}}

\newtheorem{theorem}{Theorem}
\newtheorem{proposition}{Proposition}[section]

\newtheorem{lemma}[proposition]{Lemma}
{\theorembodyfont{\rmfamily}
\newtheorem{conjecture}[proposition]{Conjecture}
\newtheorem{definition}[proposition]{Definition}

\newtheorem{example}[proposition]{Example}
\newtheorem{remark}[proposition]{Remark}
}

\include{epsf}

\begin{document}

\maketitle

{\em
\begin{flushright}
\begin{tabular}{rl}
\hfill & \\
& To Sasha Mednykh \\
& on the occasion \\
& of his 50th birthday
\end{tabular}
\end{flushright}
}

\begin{abstract}
The Hurwitz space is a compactification of the space of rational
functions of a given degree. The Lyashko-Looijenga map assigns to
a rational function the set of its critical values. It is known
that the number of ramified coverings of $\CP^1$ by $\CP^1$ with
prescribed ramification points and ramification types is related
to the degree of the Lyashko--Looijenga map on various strata
of the Hurwitz space. Here we explain how the degree of the
Lyashko-Looijenga map is related to the
intersection theory on this space. We describe the cohomology
algebra of the Hurwitz space and prove several relations between
the homology classes represented by various strata.
\end{abstract}

\section{Introduction}

\subsection{Statement of the problem}
In a series of two papers, the present one and the second
one by the second author~\cite{Zvonkine}, we continue the study of the
Hurwitz problem~\cite{Hurwitz1} concerning counting ramified coverings
of the $2$-sphere.
Roughly speaking, this problem can be formulated as follows:

{\it Given a set of fixed ramification points on the target
$2$-sphere and a set of ramification types over these points,
count the number of non-isomorphic ramified coverings $S\to S^2$
by a $2$-surface~$S$, having the prescribed ramification
types over the prescribed ramification points}.

Two coverings are considered to be isomorphic if there is
a homeomorphism of the covering surfaces taking the
first covering to the second one.
In fact, we count isomorphism classes with the weight
equal to the inverse order of the group of automorphisms
of the class; in this form the problem is more natural,
and admits a lot of applications.

In a more modern setting this problem can be reformulated
as the computation of the Gromov--Witten invariants of
the complex projective line. A detailed exposition of the history
of the problem and
a description of various approaches to its solution
can be found in~\cite{Lando}.

Hurwitz~\cite{Hurwitz2} himself gave an explicit answer
to the problem in the case, where the covering surface
is also a sphere, and there is one ramification point
with arbitrary ramification, while all others are points
of simple ramification. In spite of the importance
of the problem,
only recently Hurwitz's original results were
extended to more general cases and treatable answers were obtained:
\begin{itemize}
\item I.~Goulden and D.~Jackson~\cite{GouJac}
solved the combinatorial problem that allows one to
give the answer to Hurwitz's problem for polynomials
(that is, ramified coverings of the sphere by a sphere
having a point with a single preimage and
arbitrary ramification over other points);
later their results were reestablished in~\cite{LanZvo}
and~\cite{PanZvo} by absolutely different methods;
\item T.~Ekedahl, S.~Lando, M.~Shapiro, and A.~Vainshtein~\cite{ELSV}
gave an expression for the number of ramified coverings
of the sphere by a surface of arbitrary genus
with a single point of arbitrary ramification type
and all other points of simple ramification,
in terms of intersection indices on moduli spaces of curves
(when the covering surface is the sphere, the
answer given by the formula coincides with Hurwitz's one);
another proof of the formula was given in \cite{GraVak}.
\end{itemize}

Here we address the special case of the Hurwitz problem,
where the covering surface is also a sphere, but
in contrast to the Hurwitz case, ramification points
of arbitrary type are allowed. We are far from solving
the problem in this generality, but we suggest a general
approach, and show how it works in some special cases,
thus producing new enumerative results.

To a ramification point (= a critical value) of multiplicity
$d$ we assign a partition $\kappa$ of $d$, called the
{\em ramification type}. The elements of
$\kappa$ are the multiplicities of the critical points that
correspond to our critical value.

Recall the Hurwitz theorem.

\begin{theorem}[Hurwitz]
The number $h_{\kappa}$ of degree~$n$ coverings of the sphere by
the sphere, having the ramification type~$\kappa=(k_1,\dots,k_m)$
over one point, $d= k_1+ \dots +k_m$,
while all other ramification points are
fixed and simple, is
$$
h_\kappa=\frac{(2n-2-d)!}{(n-d-m)! \, |\Aut(\kappa)|}
\prod_{i=1}^m\frac{(k_i+1)^{k_i+1}}{(k_i+1)!}n^{n-d-3};
$$
here
$|\Aut(\kappa)|$ is the order of the automorphism group of the
partition~$\kappa$, $2n-2-d$ is the number of nondegenerate
ramification points, and $n-d-m$ is the number of simple preimages
of the multiple ramification point.
\end{theorem}

As a consequence of our approach we obtain,
in the second part of this paper~\cite{Zvonkine}, for example,
the following result.

\begin{theorem}
The number $h_{2,2}$ of degree~$n$ coverings
of the sphere by the sphere, having two double critical
points and $2n-6$ simple critical points, all the $2n-4$
critical values being fixed, is
$$
h_{2,2}=
\frac34 \, (27n^2-137n+180) \, \frac{n^{n-6}\, (2n-6)!}{(n-3)!}.
$$
\end{theorem}

As far as we know, this formula, as well as similar formulas for
some other ramification types, is new.

\subsection{The basis of the approach}

First of all we reduce the problem
to the calculation of the degree of some
map, called the Lyashko--Looijenga map.
This step is now standard. The Lyashko--Looijenga
map (below, the $LL$ map)
takes a meromorphic function to the set
of its critical values. Its source space can be chosen
in a variety of ways; here we define the $LL$ map
on the Hurwitz space $\cH_n$ constructed in~\cite{ELSV}.
It is the space of stable maps from genus zero complex curves
to the projective line, having trivial ramification over
infinity. All spaces of functions possessing
degenerate ramification are considered
as subvarieties in this space. These subvarieties
form a stratification of the Hurwitz space.

\begin{theorem}\label{Thmml}
The number $h_{\{\kappa_1,\dots,\kappa_c\}}$ of
ramified coverings of the sphere, having the
ramification types $\kappa_1,\dots,\kappa_c$
over prescribed ramification points with
multiplicities $d_1, \dots, d_c$, is given by the formula
$$
h_{\{\kappa_1,\dots,\kappa_c\}}=
\frac1{n!} \,
\frac{|\Aut \{\kappa_1, \dots, \kappa_c \}|}
{|\Aut \{d_1, \dots, d_c \}| } \,
\mu_{\{\kappa_1,\dots,\kappa_c\}}
$$
where $\mu_{\{\kappa_1,\dots,\kappa_c\}}$
is the degree of the $LL$ map restricted to the
stratum $\Sigma_{\{\kappa_1,\dots,\kappa_c\}}$
consisting of functions with these ramification types.
\end{theorem}

This is an instance of a general situation, see e.g.~\cite{Lando}.
Indeed, if we fix a ramified covering $f:S\to S^2$ and choose a
complex structure on the target sphere~$S^2$, then there exists a
unique complex structure on the covering surface~$S$ making the
function~$f$ into a meromorphic function. This complex structure
is produced by the Riemann construction. Hence, there is a
one-to-one correspondence between ramified coverings with fixed
ramification points of given types and meromorphic functions with
fixed critical values of the same types. The latter number is
exactly the degree of the $LL$ map on the
corresponding moduli space. The coefficient $1/n!$ results from
the fact that in our construction of the space $\cH_n$ we choose a
numbering of the $n$ poles of each rational function. The factor
$$
\frac{|\Aut \{\kappa_1, \dots, \kappa_c \}|}
{|\Aut \{d_1,\dots, d_c \} |}
$$
is due to the fact that if we permute the
ramification types $\kappa_1, \dots, \kappa_c$ preserving the
multiplicities $d_1, \dots, d_c$, we obtain a new set
of ramified coverings that lie on the same stratum
$\Sigma_{\{\kappa_1,\dots,\kappa_c\}}$ and have the same
image under the $LL$ map.

The space~$\cH_n$ is, in fact, a vector bundle
over the Deligne--Mumford moduli space $\ocM_{0;n}$
of stable genus zero curves. The latter space is
a smooth projective variety. The multiplicative group~$\C^*$
of nonzero complex numbers acts on this bundle
fiberwise. This action consists in just multiplying
a meromorphic function by a constant.
It preserves the stratification of the space~$\cH_n$
because the multiplication does not change the
type of singularities.
Deleting the zero section of the bundle and taking the
quotient modulo this action we reduce the variety
to the projectivization~$\P\cH_n$ of the vector bundle~$\cH_n$.

The projectivization~$\P\cH_n$
carries the tautological line bundle $\cT \to \P\cH_n$, and
a natural cohomology class $\Psi_n\in H^2(\P\cH_n)$,
$\Psi_n =c_1(\cT^\vee)$, the class of a hyperplane
section. The degree of the Lyashko--Looijenga map
restricted to a subvariety in~$\cH_n$
is related to the intersection
index of the subvariety with the complementary power
of the class $\Psi_n$. Namely, the following statement is true.

\begin{theorem} \label{Thm1.4}
The degree $\mu_{\{\kappa_1,\dots,\kappa_c\}}$
of the $LL$ map restricted to the
stratum $\Sigma_{\{\kappa_1,\dots,\kappa_c\}}$
is equal to
$$
|\Aut \{d_1, \dots, d_c \}| \,
\langle
\P\Sigma_{\{\kappa_1,\dots,\kappa_c\}},\Psi_n^d\rangle.
$$
Here~$d$ is the dimension of the stratum
$\P\Sigma_{\{\kappa_1,\dots,\kappa_c\}}$, $d_i$ is the multiplicity
of the $i$th ramification point {\rm(}or the sum of the elements of
$\kappa_i${\rm)}, and
$\left< \cdot, \cdot \right>$ is the coupling between
a homology and a cohomology class.
\end{theorem}

This theorem is proved in Section~\ref{Ssec:stratHur}.

\subsection{Cohomology of the Hurwitz spaces
and Kazarian's principle}

Theorem~\ref{Thm1.4} implies that our approach requires
the study of the cohomology ring $H^*(\P\cH_n)$.
All cohomology groups we consider are with complex coefficients,
and we do not specify the coefficients explicitly.
We need
\begin{itemize}
\item a reasonable description of this ring, say, in
terms of generators and relations;
\item reasonable expressions for the cohomology
classes of the strata and the class~$\Psi_n$
in terms of the generators.
\end{itemize}
The solution of the first problem is known
since the space~$\cH_n$ is a vector bundle
over the moduli space of stable rational curves, whose
cohomology is known. The second problem seems to be much more
difficult.

In fact, we need less than the whole cohomology
ring. The symmetric group~$S_n$
acts on the space~$\cH_n$ by permuting the indices, whence
it acts on the cohomology space~$H^*(\P\cH_n)$.
All the cohomology classes we are interested in
are symmetric with respect to this action.
Therefore, we only need to know the $S_n$-symmetric part
of the cohomology algebra~$H^*(\P\cH_n)$.
This is an additional problem, but its solution can
lead to a simplification of
the decompositions of the strata.

The situation would become even easier if we restrict
ourselves to the subalgebra in the cohomology algebra~$H^*(\P\cH_n)$
generated by the classes of the strata.

Consider two closed codimension~1 subvarieties $C_n$ and $M_n$
in $\cH_n$; the first of them is called the {\it caustic}
and it is the closure of the space of functions
having a critical point of order~$2$,
the second one is the {\it Maxwell stratum}
and it is the closure of the space of functions
having two distinct critical points with coinciding
critical values.

\begin{conjecture} \label{Conj:kazarian}
The subalgebra of the cohomology algebra $H^*(\P\cH_n)$
generated by the cohomology classes of the strata
is generated by the two classes $C_n$ and
$M_n$ in $H^2(\P\cH_n)$.
\end{conjecture}

Up to now, the basis of this conjecture is not too solid.
But it does not contradict our sample calculations, and it has
a nonformal justification coming from Kazarian's theory.

Kazarian's theory concerns cohomology classes
of multisingulatities of a map $f:M\to N$ of two complex
manifolds, $\dim M \leq \dim N$. Its main statement
claims that there is a universal way to express
these cohomology classes
in terms of the characteristic classes
of the tangent bundles over the two manifolds $M$ and $N$,
more precisely in terms of the class $c(TM)/f^*c(TN)$,
where $c$ is the total Chern class.
It is a development of the theory of the Thom polynomial.
Up to now, only a preliminary version~\cite{Kazarian}
of the text describing the theory is available,
and we are not going to refer directly to statements
from it. However, the ideology of the theory
seems to be applicable in the situation we are studying,
and it leads to a number of conjectures concerning
the part of the cohomology ring of the Hurwitz spaces
we are interested in.

In our situation, the variety $M$ is the universal curve
${\cal U}_n$ over the Hurwitz space $\P\cH_n$, while
the variety $N$ is the quotient of 
$\cH_n \times \CP^1$ by the natural $\C^*$ action. The map $f$
is the universal map over $\P\cH_n$. It is easy to
see that the map $f$ almost identifies the tangent spaces
$TM$ and $TN$; more precisely, one has $c(TM)/f^*c(TN)
= (1+a)/(1+b)$, where $a$ and $b$ are first Chern classes
of some linear bundles. However, new complications
arise, since we are interested in cohomology classes
of $\P\cH_n$ and not of ${\cal U}_n$.

Kazarian's ideology leads to Conjecture~\ref{Conj:kazarian}.

If the conjecture is true, then the computation
of intersection numbers of the strata in the Hurwitz space
with the complementary powers of the class $\Psi_n\in H^2(\P\cH_n)$
can be split into two independent stages:
\begin{itemize}
\item express the cohomology class of the required stratum
as a (homogeneous) polynomial in the classes~$M_n$ and~$C_n$;
\item find the intersection index of each monomial in~$M_n$
and~$C_n$ with the complementary degree of the class~$\Psi_n$.
\end{itemize}

Another cohomology class in $H^2(\P\cH_n)$ is also
distinguished. This is the class $\Delta_n$
dual to the subvariety of functions defined   
on singular curves. This class also can be expressed
in terms of the classes $C_n$ and $M_n$.
This means that the subring in $H^*(\P\cH_n)$
generated by the classes $C_n$ and $M_n$
coincides with that generated by
$\Psi_n$ and $\Delta_n$. The linear relations
relating these four classes are as follows:
\begin{eqnarray*}
C_n&=&6(n-1)\Psi_n-3\Delta_n; \\
M_n&=&2(n-1)(n-6)\Psi_n+4\Delta_n.
\end{eqnarray*}  

The coefficients in these linear relations are
polynomial in~$n$, which (together
with the informal support of Kazarian's theory)
allows us to sharpen Conjecture~1.5
in the following way:

\begin{conjecture}
Each cohomology class dual to a stratum
in~$\P\cH_n$ can be expressed as a homogeneous
polynomial in~$\Psi_n,\Delta_n$
whose coefficients that are polynomials in~$n$.
\end{conjecture}

Therefore the calculations must even be simplified.

\subsection{Plan of the papers}
In Section~2 we give precise definitions of notions mentioned
in this introduction. Section~3 is devoted to the
description of the Lyashko--Looijenga map.
In Section~4 the structure of cohomology of the Hurwitz
spaces is analyzed.

\subsection{Acknowledgments}

We are greatful to A.~Zorich who invited
us both to Rennes, where the main part
of this work was accomplished. The work was
completed during our participation
in the Mini-Workshop on Hurwitz theory
and its ramifications, and we thank its organizers 
J.~H.~Kwak and A.~D.~Mednykh for inviting us.   
We also thank B.~Shapiro and M.~Shapiro
for their interest and M.~Kazarian for explaining
the basics of his theory.

\section{Main definitions}

In this section we briefly recall the definitions of
stable curves and stable maps (we will restrict our attention
to curves of genus 0), as well as the definition of the
Hurwitz space~$\cH_n$ from~\cite{ELSV}.

\subsection{Moduli spaces of stable curves}

Consider the Riemann sphere $\CP^1$
with $n \geq 3$ distinct numbered marked points on it.
Two choices of marked points are
considered equivalent if there exists an automorphism
of $\CP^1$ sending one set of marked points
to the other, preserving their numbering.

The space of all
nonequivalent choices of $n$ numbered marked points
on $\CP^1$ can be endowed with a natural structure
of a smooth irreducible (noncompact) complex
manifold of dimension $n-3$.
We denote this moduli space by~$\cM_n$.

\begin{example} \label{ExM3M4}
The space $\cM_3$ is a point, because any triple of points
on $\CP^1$ can be sent to any other triple by a
biholomorphic map. Similarly, $\cM_4$ is isomorphic to
$\CP^1 - \{0,1, \infty\}$ (here we suppose that
the projective line is endowed with an arbitrary complex coordinate).
Indeed, for any four marked points
on $\CP^1$ one can send the first three to $0$,
$1$, and $\infty$. The fourth marked point will be
sent to some uniquely determined point $\lambda$,
different from $0$, $1$, and $\infty$, and this value
of $\lambda$ is the point of $\cM_4$.
\end{example}

The spaces $\cM_n$ are compactified by adding to them
new points corresponding to singular stable curves.
All the curves under consideration are compact. The only
singularities allowed are simple {\it nodes}, or {\it double points}
(as at the origin in the plane curve $xy=0$);
such curves are called {\it nodal}.

Separating the two branches of a nodal curve at each node,
we obtain a disjoint union of smooth compact
curves, which is called the {\it normalization}
of the original curve. The connected components of the
normalization are the {\it irreducible components}
of the curve. Associate to a nodal curve a
graph whose vertices correspond to the connected components
of the normalization, and two vertices are connected by an edge
iff the corresponding components intersect
at a node (loops appear if a component intersects itself).
This graph is called the {\it modular graph} of the nodal curve.
By definition, the genus of a nodal curve is zero
if all connected components of its normalization
are projective lines, and the modular graph is a tree.

\begin{definition}
A {\em stable genus 0 curve}
is a nodal curve of genus~$0$
with $n$ distinct marked and numbered points
satisfying the following conditions: (i) the marked
points do not coincide with the nodal points; (ii)
the number of automorphisms of the curve preserving the
marked points is finite.
\end{definition}

The second condition is equivalent to the following
easily verifiable condition (ii'):
the total number of marked and nodal points on every
irreducible component of the nodal curve is greater than or equal
to~$3$.

The space of genus zero stable curves with~$n$ marked
points is endowed with a natural structure of
a smooth irreducible
compact complex manifold of dimension $n-3$;
we denote this space by~$\ocM_n$.
The space~$\cM_n$ of smooth curves form a dense subvariety
in~$\ocM_n$.

\begin{example}
Apart from smooth curves, there
are exactly $3$ stable genus 0 curves with
$4$ marked points. They are shown in Fig.~\ref{FigM4}.
If $\lambda$ from Example~\ref{ExM3M4} tends to
$0$ (respectively, $1$ or $\infty$), the corresponding
smooth curve tends to the 3rd (respectively 2nd, or
1st) stable curve in the figure.

\begin{figure}
\
\begin{center}
\epsfbox{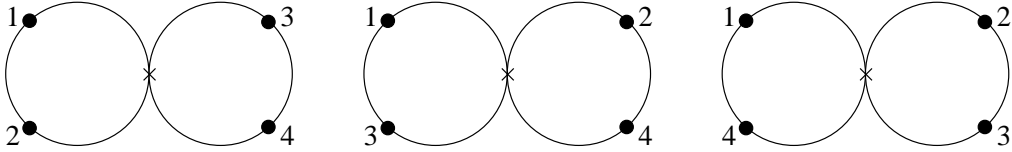}
\end{center}

\caption{Singular stable curves with $4$ marked points.
\label{FigM4}}

\end{figure}

Accordingly, $\ocM_4$ is obtained from $\cM_4$ by adding
three initially punctured points. Actually, $\ocM_4$ is isomorphic to
$\CP^1$.
\end{example}

\subsection{Hurwitz spaces}

Now consider the space of meromorphic functions on
$\CP^1$ with exactly $n \geq 1$ numbered simple poles.
This space has a structure of a smooth (noncompact)
complex manifold of dimension $2n-2$. In the case when $n \geq 3$,
if we forget the function itself
and only remember the positions of its poles, we
obtain a map from the space of functions onto $\cM_n$.
Because of this, we will often say ``marked points''
instead of ``poles of $f$''.

We are going to construct a compactification of
the space of meromorphic functions following~\cite{ELSV},
similar to that of the
moduli space $\cM_n$. To do that, we will need
to define meromorphic functions on nodal curves.
A meromorphic function $f$ on a nodal curve $S$ is simply
a meromorphic function defined on each irreducible component
of the nodal curve. If two components
intersect at a nodal point, the function must take
the same value at this point on both components.
Two pairs $(S_1,f_1)$ and $(S_2, f_2)$
are equivalent if there exists a biholomorphic
isomorphism $\varphi:S_1 \rightarrow S_2$
that preserves the numbering of the poles and
makes the following diagram commutative:

\setlength{\unitlength}{1mm}

\begin{center}

\begin{picture}(21,15)(0,0)

\put(0,12){$S_1$}
\put(16,12){$S_2$}
\put(6.5,-2){$\CP^1$}
\put(5,13.5){\vector(1,0){9}}
\put(1.5,10.5){\vector(3,-4){6}}
\put(17.5,10.5){\vector(-3,-4){6}}
\put(8,16){$\varphi$}
\put(.5,5){$f_1$}
\put(16,5){$f_2$}

\end{picture}

\end{center}

\begin{definition}
A {\em stable genus 0 meromorphic function} is
a function $f: S \rightarrow \CP^1$ defined
on a nodal curve $S$ and satisfying the
following conditions. (i) The function $f$ does not
have poles at nodal points. (ii) The number of
automorphisms of the pair $(S,f)$ is
finite.
\end{definition}

As above, a comment must be made about the
second condition. An automorphism of the pair
$(S,f)$ is a map $\varphi$ from $S$
to itself that makes the following diagram
commutative:

\setlength{\unitlength}{1mm}

\begin{center}

\begin{picture}(21,15)(0,0)

\put(0,12){$S$}
\put(16,12){$S$}
\put(6.5,-2){$\CP^1$}
\put(5,13.50){\vector(1,0){9}}
\put(1.5,10.5){\vector(3,-4){6}}
\put(17.5,10.5){\vector(-3,-4){6}}
\put(8,16){$\varphi$}
\put(.5,5){$f$}
\put(16,5){$f$}

\end{picture}

\end{center}

Condition (ii) is again easy to check,
because it is equivalent to (ii'): If the
function $f$ is constant on some irreducible component
of $S$, then the number of nodal points
on this component is greater than or equal to $3$.
(Note that such a component cannot contain marked
points precisely because $f$ is constant on it.)

\begin{definition}
The {\em Hurwitz space} $\cH_n$  (for $n \geq 2$)
is the space of all
stable genus 0 meromorphic functions with $n$
simple numbered poles.
\end{definition}

There is a natural action of $\C^*$ on $\cH_n$ defined
by simply multiplying a stable function by a constant.
The invariant points of this action are stable functions
such that the nodal curve has at most one pole on every component
and the function vanishes on the components without poles.

These points form the {\em zero section} of the Hurwitz
space. (We will soon see that the Hurwitz space is
a vector bundle and the points in question indeed form
its zero section.)

\begin{definition}
The {\em projectivized Hurwitz space} $\P\cH_n$ is
the quotient of $\cH_n$, without the zero section,
by $\C^*$.
\end{definition}

For $n \geq 3$, $\P\cH_n$ has the structure of a smooth
compact complex manifold of dimension $2n-3$.
$\P\cH_2$ is a compact smooth orbifold.

\begin{example}
$\P\cH_2$ is isomorphic (as an orbifold) to the weighted projective
line with weights $2$ and $1$, that is to
$$
(\C^2 - \{ 0 \})/ (z,w) \sim (\lambda^2 z, \lambda w).
$$
Indeed, if we are given a meromorphic map
with 2 poles on $\CP^1$, we can move its poles to
$0$ and $\infty$ and write the function as $f(z) = az+b+c/z$.
The numbers $a,b,c$ are considered up to (simultaneous)
multiplication by a scalar factor.
Moreover, making a change of variables $z \mapsto \lambda z$,
we see that the function $(\lambda a) z + b + (c/\lambda)/z$
is equivalent to $f$. Thus $[ac:b]$ is a set of homogeneous
coordinates in $\P\cH_2$. The point $ac=0$
corresponds to the stable meromorphic function
defined on a nodal curve with two irreducible components.

The space $\P\cH_3$, as will soon be proved, is isomorphic to
$\CP^3$.
\end{example}

\subsection{$\cH_n$ is a bundle over $\ocM_n$} \label{Ssec:bundle}

Now we are going to prove that for $n \geq 3$, $\P\cH_n$ is
the projectivization of an $(n+1)$-dimensional
holomorphic vector bundle on $\ocM_n$. First
let us define several line bundles on $\ocM_n$.

Consider a stable curve $S$ with $n$ marked
points. For $i$ from $1$ to $n$,
let $L_i$ and $L_i^\vee$ be the complex line respectively
cotangent and tangent to $S$ at the $i$th marked point.

\begin{definition}\label{DefL_i}
$\L_i$ (respectively $\L^\vee_i$) is the line bundle over $\ocM_n$
whose fiber over a point $S \in \ocM_n$ is the
line $L_i$ (respectively $L_i^\vee$).
\end{definition}

Denote by $\C$ the trivial line bundle and
let $E_n$ momentarily denote the following
$(n+1)$-dimensional bundle over $\ocM_n$:
$$
E_n = \L_1^\vee + \dots + \L_n^\vee + \C.
$$

\begin{proposition}\label{PropH_n=P(E_n)}
For $n \geq 3$, we have $\cH_n = E_n$.
\end{proposition}

\paragraph{Proof.} We will show that giving a stable
meromorphic function on a nodal curve is equivalent
to giving a point in the total space of the bundle
$E_n$. Considering stable functions up to a scalar
factor corresponds to taking the projectivization of $E_n$.

First of all, note that if a meromorphic function $f$ has
a simple pole at some smooth point $z$ of a nodal
curve $S$, then it determines a tangent vector to $S$
at $z$; more precisely, the principal part of the pole
is such a tangent vector. Indeed, if we want to
couple this tangent vector with a cotangent
vector represented by a $1$-form $\omega$,
the result will be the residue of $f\omega$ at $z$.

Now, given a stable meromorphic function on a nodal
curve, we can divide its poles into two parts. The poles
of the first type lie on the irreducible components of the curve
that contain at least $3$ marked and double points.
The poles of the second type lie on the irreducible components
that contain a unique pole and a unique nodal point.
(These are the only possible cases.)

If we forget the function itself
and retain only the positions of its poles, we
obtain a nodal curve with $n$ marked points. This
curve is not necessarily stable. However, we can
make it stable by contracting every component
that contains only 1 nodal and 1 marked point.
The former node becomes the marked point.
Thus we obtain a map from $\cH_n$ to $\ocM_n$.

Further, consider the stable curve that we have obtained. We
assign to its marked points of the first type the principal parts
of the corresponding poles (they are tangent vectors at the marked
points). We assign zero tangent vectors to the marked points of
the second type. Finally, we assign to the stable meromorphic
function a complex number: the sum of its critical values given by
the Lyashko-Looijenga map (see Section~\ref{Ssec:LL} below).
Thus we have assigned to the stable meromorphic function a point
in the total space of $E_n$.

Conversely, given a point in the total space of $E_n$,
we can easily recover the stable meromorphic function.
First we extend the stable curve by adding new irreducible
components at each marked point of the second type
(those, whose tangent vector is equal to 0). Then,
knowing the principal parts of the poles, we can
recover the function $f$ uniquely up to an additive
constant. Knowing the sum of its critical values
allows us to fix the constant. \qed

\begin{remark}
The vector bundle $\cH_n$ can be understood more 
conceptually. 

Consider the space ${\cal S}_n(\CP^1)$ of all stable
degree $n$ maps from genus $0$ curves $S$ to $\CP^1$. Inside
this space there is a suborbifold canonically isomorphic
to $\ocM_n/S_n$, where $S_n$ is the symmetric group acting
on $\ocM_n$ by renumbering the poles. This subvariety
consists of stable maps of the following form. Take a
stable genus $0$ curve $S'$ with $n$ marked but not numbered
points. Attach to $S'$ a new spherical component at each 
marked point. Instead of the former marked points, mark
one point on each of the new spherical components.
Thus we have obtained a nodal curve $S$. Consider the
map $f:S \rightarrow \CP^1$ that sends the whole curve
$S'$ to $0$ and is of degree one on each of the new
spherical components, the poles being at the new marked
points. Such stable maps $(S,f)$ form an orbifold
isomorphic to $\ocM_n/S_n$. 

Now, the normal bundle to $\ocM_n/S_n$ in ${\cal S}_n(\CP^1)$ 
(in the orbifold sense) is canonically identified with
the quotient $\cH_n/S_n$.
\end{remark}

\begin{definition}\label{Def:Psi}
We denote by $\cT$ the canonical line bundle over $\P\cH_n$
and by $\Psi_n = c_1 (\cT^\vee) \in H^2(\P\cH_n)$ 
the first Chern class of the dual line bundle.
\end{definition}

\section{The Lyashko-Looijenga map} \label{Ssec:LL}

\subsection{The definition of the map}

Consider the space of unordered sets of $m$ (not
necessarily distinct) complex
numbers. This space is isomorphic to the space
of monic polynomials of degree~$m$: just take
a set of~$m$ numbers to the polynomial whose
roots are these numbers.

\begin{definition}
The {\em Lyashko-Looijenga map} $LL$
is the map that assigns to a meromorphic function with
$n$ simple poles on $\CP^1$ the (unordered) set of
its $2n-2$ critical values counted with multiplicities.
\end{definition}

We recall that a {\em critical point} of a function $f$
is a point, where $df=0$, while a {\em
critical value} of $f$ is its value at a critical point.
The {\em multiplicity of a critical point} is the
multiplicity of the zero of $df$ at this point. In
other words, if the function $f$ has the form $f(z) = z^k + c$
in some local coordinate $z$, then $z=0$ is a critical
point of multiplicity $k-1$. The
{\em multiplicity of a critical value} is the sum
of multiplicities of the critical points at which the
critical value in question is attained. If we consider
the function as a ramified covering of $\CP^1$,
then its critical values are exactly those
points in the image over which the covering is ramified.

The notion of a {\it cone} (or of a cone bundle) generalizes
that of a vector bundle to the case
of singular fibers, see~\cite{Fulton}. The essential point is
that cones carry a fiberwise action of the multiplicative group~$\C^*$
of nonzero complex numbers. A cone morphism preserves this action.

\begin{proposition}[\cite{ELSV}]
The Lyashko-Looijenga map extends to a cone morphism
$$
LL:\cH_n \rightarrow \C^{2n-2}.
$$
\end{proposition}

The proof of the proposition is based on
the extension of the notion of set of critical values to
maps defined on nodal curves.

First consider all the irreducible
components of $S$ where $f$ is constant. Consider the
union of these irreducible components, and let $S_0$ be
a connected component of this union. (In general
$S_0$ consists of several irreducible components
intersecting at double points.) Suppose $S_0$ contains
$k$ nodal points at which it intersects other
irreducible components (on which $f$ is not constant).
Then the value of $f$ on $S_0$ is considered a
critical value of multiplicity $2k-2$.

Second, consider a nodal point lying on two
irreducible components of $S$ such that
$f$ is not constant on either of them. The
value that $f$ takes at such a point is considered
a double critical value.

Finally, the critical values of $f$ on all the
irreducible components are taken into account
as ordinary critical values.

One can prove that in that way we have obtained
a set of $2n-2$ critical values, and that this
set depends continuously on the point of $\cH_n$.
See~\cite{ELSV} for more details.

\subsection{A stratification of the Hurwitz space}
\label{Ssec:stratHur}

Consider a degree~$n$ ramified covering $f:S\to S^2$
of the sphere. Let $t\in S^2$
be a ramification point of $f$ of multiplicity $d$. To
this ramification point we can assign
a partition~$\kappa$ of~$d$.
This partition consists of the multiplicities of all
critical points of $f$ in $S$ whose image under $f$ is equal to $t$.
We will write such a partition in the multiplicative form,
$\kappa=1^{m_1}\dots d^{m_d}$, where $m_i$,
$i=1,\dots,d$ is the number of critical points in the
preimage of $t$ having multiplicity~$i$,
$1\cdot m_1+\dots+d\cdot m_d = d(\kappa) = d$.
If the covering has~$c$ ramification points, $t_1,\dots,t_c$,
then~$c$ partitions~$\kappa_1,\dots,\kappa_c$
are associated to it.

Sets of partitions of~$n$ that can appear
as sets of ramification partitions of a ramified covering
are subject to a number of constraints.
If the covering surface~$S$ also is the sphere,
the Riemann--Hurwitz formula implies that the total
degeneracy of the set of ramification partitions
of a meromorphic function must be~$2n-2$,
$$
d(\kappa_1)+\dots+d(\kappa_c)=2n-2.
$$
However, this property alone does not guarantee
the existence of a ramified covering with
the prescribed ramification type.

\begin{definition}\label{Defstrata}
The set of all rational functions on $\CP^1$ with $n$
simple numbered poles such that the set of partitions
assigned to their critical values is
$\{ \kappa_1, \dots, \kappa_c \}$ is called an {\em open stratum}.
Its closure in $\cH_n$ is called just a {\em stratum} and denoted by
$\Sigma_{\{ \kappa_1, \dots, \kappa_c \}}$.
\end{definition}

The projectivizations of these strata (lying in $\P\cH_n$) are
defined in the obvious way. They will be denoted
by $\P \Sigma_{\{ \kappa_1, \dots, \kappa_c \}}$.

Sometimes we will omit reference to simple partitions $\kappa = 1$
in the index, omitting in this case the braces as well.

In the sequel we will often describe strata by specifying
the set of functions that belong to their open parts and omitting
the words ``closure of''. The reader should bear in mind
that we always consider closed strata, unless explicitly
stated otherwise.

\begin{example} \label{ExCM}
There are two strata of complex codimension~1. They are
the caustic $C_n\subset\cH_n$ and
the Maxwell stratum $M_n\subset\cH_n$. The
{\em caustic} $C_n = \Sigma_{2}$ is (the closure of) the set of functions
with a double critical point. The {\em Maxwell stratum}
$M_n = \Sigma_{1^2}$ is (the closure of) the set of functions
with the same critical value taken at two distinct critical
points.
\end{example}

\begin{definition}\label{Defmultiplicity}
The {\em degree} of the $LL$ map on a stratum
$\Sigma_{\{ \kappa_1, \dots, \kappa_c \}}$
is the number of the preimages in
$\Sigma_{\{ \kappa_1, \dots, \kappa_c \}}$ of a
generic point in $LL(\Sigma_{\{ \kappa_1, \dots, \kappa_c \}})$.
The degree is
denoted by $\mu_{\{ \kappa_1, \dots, \kappa_c \}}$.
\end{definition}

\begin{example}
Since the space $\ocM_3$ is a point, the space $\cH_3$ is
isomorphic to $\C^4 = L_1^\vee + L_2^\vee + L_3^\vee + \C$,
where the $L_i$ are the cotangent lines at the marked
points on $\CP^1$. It can be identified with the
space of rational functions of the form
$$
az-\frac{b}z-\frac{c}{z-1} + d.
$$

In Figure~\ref{H3} we represent the stratification
of the projectivized Hurwitz space~$\P\cH_3$. The figure
shows the real part of this space with a slight simplification:
Instead of a 3-dimensional projective space we have
represented a 2-dimensional projective space, dropping
the parameter $d$.

The interesting subvarieties
in $\P\cH_3$ are the following ones. First, the
projectivized caustic~$\P C_3$, which is
a singular cubic with one self-intersection point.
It has the homogeneous equation
$$
a^3+b^3+c^3+3ab^2+3ac^2+3ba^2+3bc^2+3ca^2+3cb^2-21abc=0.
$$
Since our picture shows the real
part of $\P\cH_3$, the self-intersection point,
which has the homogeneous coordinates~$(1:1:1)$, appears to lie
apart from the rest of the cubic. Second, the three projective
lines $a=0$, $b=0$, $c=0$,
which correspond to reducible curves with one pole
lying on one component and the other two poles on the
other one. The Maxwell stratum is empty. The cubic
intersects each of the three lines at one point, but
the intersection has multiplicity three. These three
points lie on the line $a+b+c=0$.

The boxes around the picture of $\P\cH_3$ show, for each
of the subvarieties, what are the corresponding stable
rational functions. Short traits represent poles,
crosses represent critical points, while crosses
with double lines represent double critical points.
The numbers written near the arrows are the multiplicities
of the $LL$ map on the corresponding subvariety.

\begin{figure}
\begin{center}
\
\epsfbox{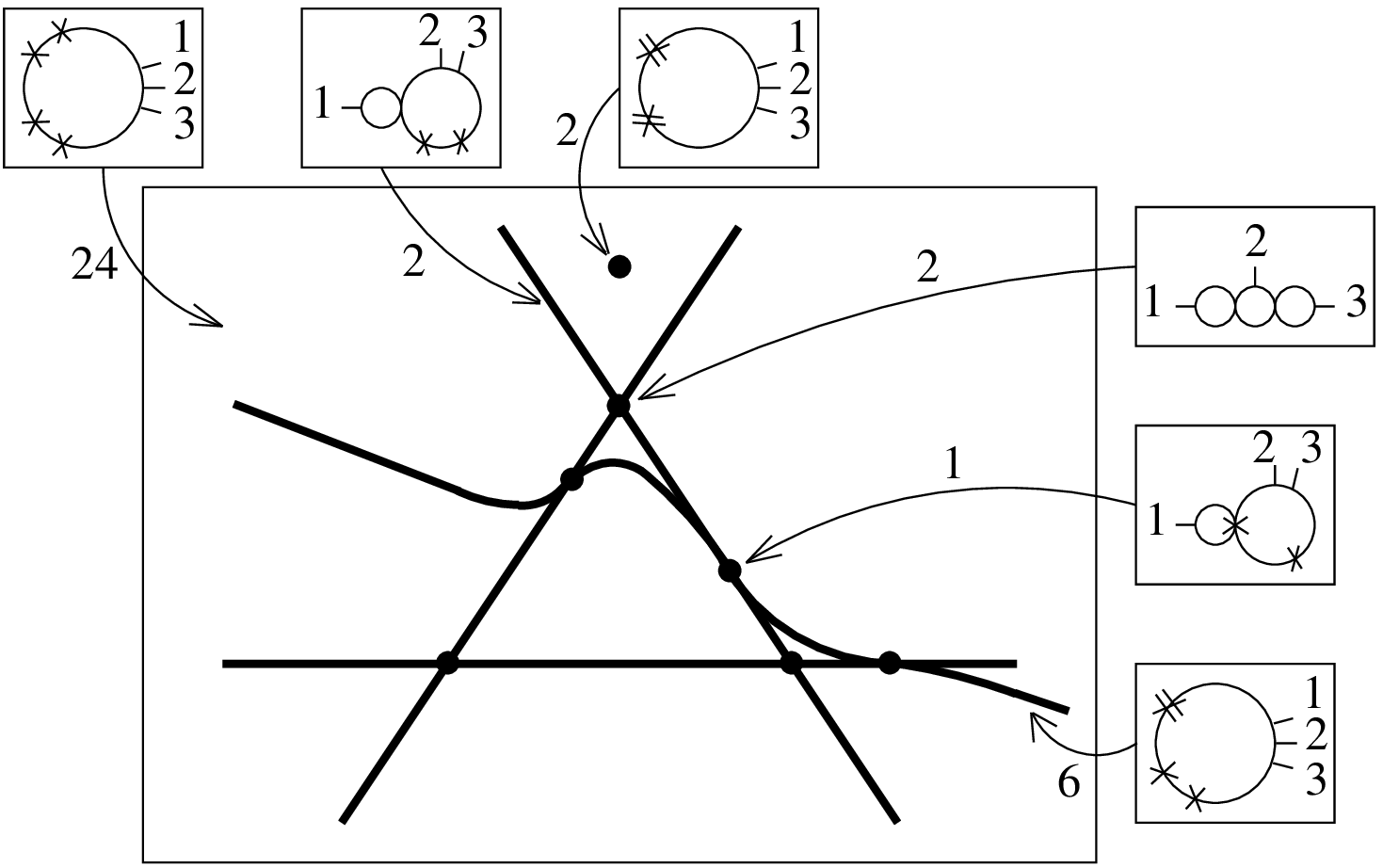}

\caption{The stratification of the ``simplified'' projectivized
space $\P\cH_3$. \label{H3}}
\end{center}
\end{figure}
\end{example}

Now we are able to prove Theorem~\ref{Thm1.4} from the
introduction. We restate it here.

Consider a stratum $\Sigma_{\{\kappa_1,\dots,\kappa_c\}}\subset \cH_n$
(see Definitions \ref{Defstrata}
and~\ref{Defmultiplicity}) determined by partitions
$\kappa_1, \dots, \kappa_c$. Let $d_1 = d(\kappa_1), \dots,
d_c = d(\kappa_c)$ be the multiplicities of the critical
values. Let
$|\Aut \{d_1,\dots,d_c\}|$ be the number
of automorphisms of the set $\{d_1,\dots, d_c\}$, i.e.,
the number of permutations $\sigma$ of $c$ elements such that
$d_i = d_{\sigma(i)}$ for all $i$. 

\setcounter{theorem}{3}

\begin{theorem}\label{Thmmult}
The degree $\mu_{\{\kappa_1,\dots,\kappa_c\}}$
of the Lyashko-Looijenga map on
$\Sigma_{\{\kappa_1,\dots,\kappa_c\}}$ equals
$$
|\Aut \{d_1, \dots, d_c \}| \,
\langle
\P\Sigma_{\{\kappa_1,\dots,\kappa_c\}},\Psi_n^d\rangle,
$$
where~$d$ is the dimension of the stratum
$\P\Sigma_{\{\kappa_1,\dots,\kappa_c\}}$ and
$\left< \cdot, \cdot \right>$ is the coupling between
a homology and a cohomology class.
\end{theorem}

\paragraph{Proof of Theorem~\ref{Thmmult}.}
To prove the theorem we consider a $(2n-2)!$-sheeted ramified covering
$\tPcH_n$ of $\P\cH_n$. By definition, a point of $\tPcH_n$ corresponds
to a stable meromorphic function {\em with numbered critical
points}.

More precisely, consider the Zariski open subset $X$
of $\P\cH_n$ constituted by  functions $f$
defined on $\CP^1$, with $2n-2$ simple critical points.
Over this open subset we can construct a $(2n-2)!$-sheeted
non-ramified covering $\widetilde X$, whose points correspond to
functions with numbered critical points. It is a connected
complex manifold. Then this covering is uniquely
extended to a ramified covering of $\P\cH_n$ in the following way.

For simplicity we start by describing the set of
preimages of every point $f_0 \in \P\cH_n$ in the covering
$\tPcH_n$.
Consider an $f_0 \in \P\cH_n$ and a small ball $B$
surrounding $f_0$. A point $f \in B \cap X$ has
$(2n-2)!$ preimages in the covering. The local monodromy group
$\pi_1(X \cap B)$ acts by permutations on these preimages.
The preimages of $f_0$ in the covering
are the orbits of this action. Now let us give
an algebro-geometric definition.

Consider an affine chart $U$ in $\P\cH_n$ and let $A$ be the
algebra of algebraic functions on $U$, $U = \mbox{spec }A$.
Let $F \in A$,
be a function whose set of zeroes is the ramification divisor.
Then the localization $A_F$ is the algebra of algebraic functions
on the complement $U \cap X$ to the ramification divisor. The map
$A \hookrightarrow A_F$ is an injection, because $A$ has no zero
divisors. Let $B$ be the algebra of algebraic functions on the
lifting $\widetilde {U \cap X}$ of $U \cap X$ on $\widetilde X$.
Then $A_F \hookrightarrow B$ is again an injection. The composition
of the two injections gives an injection $A \hookrightarrow B$.
We denote by $C$ the {\em integral closure} of $A$ in $B$. By definition,
the ramified covering over $U$ is then the variety $\mbox{Spec }C$. Such
ramified coverings of affine charts can then be glued into a covering
of the total space $\P\cH_n$. 

\begin{lemma} \label{Lem:normal}
The variety $\tPcH_n$ is normal.
\end{lemma}

\paragraph{Proof.} Recall that a normal algebraic variety
is a variety such that each of its points
possesses an affine neighborhood whose algebra of functions
is integrally closed in its field of fractions
(see~\cite{Greco}). On normal
varieties Chern classes of vector bundles and local
intersection indices of subvarieties are well-defined.
Normality is a local property, therefore we can restrict
our attention to an affine chart $U$ as above.

Let us prove that $C$ is integrally closed. First of all,
$B$ is integrally closed, because it is an algebra of functions
on a smooth variety. Thus an element $x$ of the field of fractions
of $C$ integral over $C$ belongs to $B$. But $C$ is integrally
closed in $B$ (being the integral closure of $A$), thus $x$
necessarily belongs to $C$.
\qed

In general, $\tPcH_n$ is not smooth. For example,
the $2$-sheeted ramified covering of $\C^2$ ramified
over the axes $x=0$ and $y=0$ is the cone given by the
equation $z^2 = xy$.

Our description of the preimage of a point in
$\P\cH_n$ under the projection $\tPcH_n \rightarrow \P\cH_n$
follows from the fact
that a normal complex algebraic variety is locally irreducible
as a complex manifold (\cite{Greco}, Theorem~6.6).

We denote by $\tLL$ the map
$$
\tLL: \tPcH_n \rightarrow \CP^{2n-3}
$$
that assigns to a meromorphic function with numbered
critical points in $\tPcH_n$ the set of its numbered critical
values. (Recall that both the function and the set of its
critical values are considered up to a common scalar factor.)
Denote by $\cO(-1)$ the tautological line bundle over
the target space $\CP^{2n-3}$.

Recall that $\cT$ is the tautological line bundle over $\P\cH_n$.
Denote by $\tcT$ the pull-back of $\cT$
under the projection $\tPcH_n\to\P\cH_n$.
Then $\tcT$ coincides with
the pull-back of $\cO(-1)$ under $\tLL$. Indeed, consider
a point of the total space of the bundle $\tcT$,
outside its zero section.  It corresponds to
a nonvanishing meromorphic function with numbered critical points
(and, this time, not considered up to a scalar factor).
Similarly, the total space of the bundle $\cO(-1)$ without
the zero section is $\C^{2n-2} - \{0\}$, and it corresponds
to sets of $2n-2$ numbered critical values, not all of which
are equal to $0$. It follows that the map $\tLL$ can be
lifted to the total spaces of the bundles $\tcT$
and $\cO(-1)$: it suffices to assign to a meromorphic function the
set of its critical values. This lifting identifies the
fibers of the two bundles.

It follows that for any
subvariety $\tPSigma \subset \tPcH_n$ of pure dimension $d$
such that its image $\tLL(\tSigma)$ is irreducible and
also of dimension $d$, we have
$$
\left< \tPSigma, \, (c_1(\tcT^\vee))^d \right>
= \mu(\tLL|_{\tPSigma}) \deg \tLL(\tPSigma).
$$
Indeed, we have
$$
\deg \tLL (\tPSigma) = \tLL(\tPSigma) \cap P,
$$
where $P\subset\CP^{2n-3}$
is a generic projective subspace of dimension
complementary to that of $\tLL(\tPSigma)$. Now, $P$ is the intersection
of $d$ projective hyperplanes. A hyperplane is a
section of $\cO (1)$ (the line bundle dual to the
tautological line bundle $\cO(-1)$). Let us count
in two different ways the number of preimages in $\tPSigma$
of the points of intersection between $\tLL(\tPSigma)$ and $P$.
First, since each intersection point has
$\mu(\tLL|_{\tPSigma})$ preimages, the number of preimages
equals
$$
\mu(\tLL|_{\tPSigma}) \deg \tLL(\tPSigma).
$$
Second, the pull-back to $\tPcH_n$ of a generic hyperplane
is a generic section of the line bundle $\tcT^\vee$ dual to $\tcT$.
Taking the intersection between $d$ sections like that
and the subvariety $\tPSigma$, we obtain
$$
\left< \tPSigma, \, (c_1(\tcT^\vee))^d \right>.
$$
Hence, these two numbers are equal.

The same equality holds if $\tLL(\tPSigma)$ is not irreducible, but
all generic points of $\tLL(\tPSigma)$ have the same number of
preimages.

Now let $\tPSigma$ be the lifting to $\tPcH_n$
of a stratum $\P\Sigma_{\{\kappa_1,\dots,\kappa_c\}}$
in $\P\cH_n$. Denote by $d$ their common dimension.
For a partition~$\kappa$
denote by $\kappa!$ the product of the factorials
of its elements.
It is easy to see that a generic point in
$\P\Sigma_{\{\kappa_1,\dots,\kappa_c\}}$ has
$$
\frac{(2n-2)!}{\kappa_1!\dots\kappa_c!}
$$
liftings to $\tSigma$. Further, the image of $\tPSigma$ under
the lifted Lyashko-Looijenga map $\tLL$ is a union
of
$$
\frac{(2n-2)!}{d_1!\dots d_c! \, |\Aut \{d_1,\dots,d_c\}|}
$$
projective subspaces of dimension $d$.

Finally, denote by
$$
\tilde \mu_{\{\kappa_1,\dots,\kappa_c\}} =
\mu(\tLL|_{\tPSigma_{\{\kappa_1,\dots,\kappa_c\}}})
$$
the degree of $\tLL$ on $\tPSigma_{\{\kappa_1,\dots,\kappa_c\}}$,
which is the number of preimages
in $\tPSigma_{\{\kappa_1,\dots,\kappa_c\}}$ of a generic point in the image of
$\tPSigma_{\{\kappa_1,\dots,\kappa_c\}}$.
Simple combinatorial considerations show that
$$
\tilde \mu_{\{\kappa_1,\dots,\kappa_c\}}
=
\prod_{i=1}^c \frac{d_i!}{\kappa_i!} \, \cdot \,
\mu_{\{\kappa_1,\dots,\kappa_c\}}.
$$
The factor $\prod \frac{d_i!}{\kappa_i!}$ is the number of ways
to number the critical points once the critical values are
numbered.

Putting everything together we see that
$$
\mu_{\{\kappa_1,\dots,\kappa_c\}}
=
\prod_{i=1}^c \frac{\kappa_i!}{d_i!} \, \cdot \,
\tilde \mu_{\{\kappa_1,\dots,\kappa_c\}}
$$
$$
= \frac{\prod \kappa_i!}{\prod d_i!} \, \cdot \,
\frac{\prod d_i! \, |\Aut \{ d_1, \dots, d_c\}|}{(2n-2)!}
\, \cdot \,
\bigl< \tPSigma , (c_1(\tcT^\vee))^d \bigr>
$$
$$
= |\Aut\{d_1,\dots,d_c\}| \,\,\,
\left<\P\Sigma_{\{\kappa_1,\dots,\kappa_c\}} \, , \,
(c_1(\cT^\vee))^d \right>
$$
$$
= |\Aut\{d_1,\dots,d_c\}| \,\,\,
\left< \P\Sigma_{\{\kappa_1,\dots,\kappa_c\}} \, , \,
\Psi_n^d \right>.
$$
This proves the theorem. \qed

\section{The cohomology of the Hurwitz space}

In this section we describe the cohomology algebra
of the moduli spaces and of the Hurwitz spaces.
We also prove several cohomological identities
which will  be used in the subsequent paper~\cite{Zvonkine}
to derive recurrence
relations on Hurwitz numbers.

\subsection{Keel's description of $H^*(\ocM_n)$}
\label{Ssec:Keel}

The description of the cohomology algebra
of the space $\ocM_n$ given in this section is
due to Keel~\cite{Keel}.

\begin{proposition}[\cite{Keel}]
The cohomology ring $H^*(\ocM_n)$ is generated
by $H^2(\ocM_n)$.
\end{proposition}

Now we will describe a set of $2$-cohomology classes
that span $H^2(\ocM_n)$. (More precisely, we will
describe cycles that represent
their Poincar\'e dual homology classes.)

Denote by $D = A \sqcup B$
an unordered partition of the set $\{1, \dots, n \}$
of marked points into two disjoint parts (a partition $A \sqcup B$
coincides with $B \sqcup A$). Each part must
contain at least $2$ points.

Consider all stable curves with exactly two irreducible
components such that one component contains
the marked points of the set $A$ and the other one
of the set $B$. Such stable curves form a (complex)
codimension~1 subvariety of $\ocM_n$.
Its closure is a smooth closed subvariety of $\ocM_n$
of complex codimension~1. Therefore it represents a
2-cohomology class which we will denote by $[D]$.

\begin{proposition}[\cite{Keel}]
The classes $[D]$ span $H^2(\ocM_n)$.
\end{proposition}

The generators $[D]$ are not linearly independent.
For example, for $n=4$ all the three generators
corresponding to the partitions
$$
D_1 = \{1,2\} \sqcup \{3,4\} \,\, , \quad
D_2 = \{1,3\} \sqcup \{2,4\} \,\, , \quad
D_3 = \{1,4\} \sqcup \{2,3\}
$$
represent the same cohomology class, (Poincar\'e dual to)
a point.

More generally, let us fix four distinct numbers $i,j,k,l$
between $1$ and $n$. Denote by $[ij\D kl]$ the sum
of all the generators $[D]\in H^2(\ocM_n)$
corresponding to partitions, where
$i$ and $j$ belong
to one of the two parts, while
$k$ and $l$ belong to the other part. Then the generators
$[D]$ satisfy the relations
$$
R_{ijkl}: \quad [ij\D kl] = [ik\D jl] = [il\D jk].
$$
These relations can be deduced using the forgetful
map $\ocM_n \rightarrow \ocM_4$ which forgets all the
marked points except for $i,j,k,l$ and contracts all
the components of the curve that have become unstable.

\begin{proposition}[\cite{Keel}]
The relations $R_{ijkl}$ span all the linear relations
on the generators $[D]$.
\end{proposition}

Finally, we say that two partitions $D_1$ and $D_2$
are {\em compatible} if the set $\{1,\dots, n\}$
can be divided into a disjoint union of three sets
$A$, $B$, $C$, in such a way that
$$
D_1 = (A \sqcup B) \sqcup C, \quad
D_2 = A \sqcup (B \sqcup C).
$$

It is easy to see that if two partitions $D_1$ and $D_2$
are not compatible, then the geometric intersection
of the corresponding cycles has codimension greater
than $2$. Therefore, the intersection of two such generators
in the cohomology algebra vanishes.

\begin{theorem}[\cite{Keel}]
The cohomology algebra $H^*(\ocM_n)$ is generated by the
generators $[D]$ modulo two kinds of relations: {\rm (i)}
the linear relations $R_{ijkl}$, {\rm (ii)} the multiplicative
relations $[D_1] [D_2] = 0$ for incompatible partitions
$D_1$ and $D_2$.
\end{theorem}

This theorem gives a complete description of the cohomology
algebra of $\ocM_n$, although in practice it leads to rather
heavy computations.

Using this theorem, Kontsevich and Manin~\cite{KonMan}
found linear generators and relations of
$H^k(\ocM_n)$ for all $k$.

\subsection{A description of the cohomology of Hurwitz spaces}
\label{Ssec:H*(Hn)}

In what follows, we will require
particular $2$-cohomology classes
$\psi_i\in H^2(\ocM_n)$.
Recall that~$\cL_i$ denotes the line bundle over~$\ocM_n$
whose fiber coincides with the cotangent line at the $i$th
marked point.

\begin{definition}\label{Defpsi_i}
We denote by $\psi_i = c_1(\cL_i)$ the first
Chern class of the line bundle $\cL_i$.
\end{definition}

For $n \geq 6$ the classes $\psi_i$ {\em do not} span $H^2(\ocM_n)$.
The following proposition expresses them in terms of the
classes $[D]$.

Let $i,j,k \in \{1, \dots, n\}$. We denote
by $[i*\D jk]$ the sum of all the classes $[D]$
corresponding to partitions containing~$i$ in one part,
and~$j,k$ in the other part. The star in the notation
reminds that an irreducible component
with a single node cannot contain a single
marked point.

\begin{proposition}\label{Proppsi_i}
For each~$i$ and any distinct $j,k$
different from~$i$ we have $\psi_i = [i*\D jk].$
\end{proposition}

\paragraph{Proof.} We construct a holomorphic section of the line
bundle $\cL_i$ in the following way. Consider a stable curve
such that the marked points $i,j,k$ are on the same
irreducible component. On this component there exists
a unique meromorphic differential with simple poles
at the points $j$ and $k$, whose residues at these poles are equal
to $1$ and $-1$ respectively.
The value of this differential at the point
$i$ is an element of the line $L_i$. Thus
we obtain a section of the bundle $\cL_i$, and it
extends in a unique way to a
meromorphic section over the entire base~$\ocM_n$.
It is easy to see that this section has no poles
(thus the section is, actually, holomorphic). Its zeroes
are precisely the classes that add up to $i*\D jk$,
and it intersects the zero section of~$\cL_i$
transversally.
\qed

\bigskip

Since~$\cH_n$ is a vector bundle over~$\ocM_n$,
the cohomology algebra of $\ocM_n$ can be
canonically seen as a subalgebra of the cohomology
algebra of $\P\cH_n$. We will therefore use the same
notation for cohomology classes on $\ocM_n$ and their
pull-backs on $\P\cH_n$.

Now we express the cohomology algebra
of $\P\cH_n$ in terms of that of $\ocM_n$.

\begin{theorem}\label{Thm:H*(H)}
$$
H^*(\P\cH_n)=H^*(\ocM_n)[\Psi_n]\,\,  / \,\,
(\Psi_n-\psi_1) \dots (\Psi_n-\psi_n) \Psi_n.
$$
\end{theorem}

\paragraph{Proof.} The assertion follows immediately
from the description of $\cH_n$ as a vector bundle
over $\ocM_n$ (Section~\ref{Ssec:bundle}) and from the following
well-known fact (see, for example,~\cite{Fulton}).
For any rank~$n$ complex vector bundle~$E$ over a complex
manifold~$M$ we have
$$
H^*(\P E) =
H^*(M) [\Psi]\,\,  / \,\,
\Bigl(\Psi^{n} + c_1(E) \Psi^{n-1} + \dots + c_{n}(E) \Bigr),
$$
that is, the cohomology algebra of~$\P E$
is isomorphic to the algebra of polynomials in one variable~$\Psi$
with coefficients in~$H^*(M)$, modulo the Chern polynomial
in~$\Psi$. Here the new variable~$\Psi$ can be identified
with the first Chern class of the tautological line bundle
over~$\P E$. \qed

This gives an explicit description
of the cohomology algebra of $\P\cH_n$ for $n \geq 3$.

\bigskip

Now, let $D = A \sqcup B$ be a partition of the set
$\{1, \dots, n \}$ into a disjoint union of two nonempty
subsets. The subsets are no longer required to have at
least 2 elements. Consider the set of stable
meromorphic functions defined on nodal curves with
$2$ irreducible components, the first component containing the
marked points of the subset $A$ and the second one
of the subset $B$. This set is an open submanifold
of $\cH_n$ of complex codimension~$1$, so its closure
is a cycle Poincar\'e dual to a cohomology class
in $H^2(\P\cH_n)$. This class will be denoted by $[D]$.

If $A$ and $B$ both contain at least two elements,
then the corresponding class $[D]$ is a pull-back from
a $2$-cohomology class on $\ocM_n$.

Let $i,j,k \in \{1,\dots,n\}$.
Denote by $[i\D jk]\in H^2(\P\cH_n)$ the sum of all the classes $[D]$
corresponding to partitions
such that $i$ belongs to one part and $j$ and $k$ belong to
the other one. (There is no star after $i$, because
the singleton $\{ i \}$ is now an allowed subset.)

\begin{proposition} \label{PropiDjk}
For any pairwise distinct $i,j,k \in \{1,\dots,n\}$,
we have $\Psi_n = [i\D jk]$.
\end{proposition}

\paragraph{Proof.}
Recall that $\Psi_n$
is the first Chern class of the
tautological line bundle on $\P\cH_n$.
A section of the vector bundle $\cH_n^\vee$ dual
to $\cH_n$ will automatically provide a section of
$\cO_s(-1)$ (because a linear form on each fiber of
$\cH_n$ can be restricted to any line contained in that
fiber). But a section of the bundle $\cH_n^\vee$ can be
obtained from a section of any of the line bundles $\cL_i$
by taking its direct sum with zero sections
of each of the other line bundles~$\cL_j$
for $j\ne i$.
Thus the class $\Psi_n$ is equal to the
sum of the class $\psi_i = c_1(\L_i)$ and of the
class Poincar\'e dual to the projectivization
of the subbundle
$$
\L_1^\vee + \dots + \L_{i-1}^\vee + \L_{i+1}^\vee
+ \dots + \L_n^\vee + \C.
$$
The latter $2$-cohomology class is exactly the class
that corresponds to the partition
$$
D = \{i\} \sqcup \{1, \dots, i-1, i+1, \dots n\}.
$$
Using the expression for $\psi_i$ given in
Proposition~\ref{Proppsi_i} we obtain
$\Psi_n = [i\D jk]$.
\qed

\begin{proposition}
For $n \geq 2$, the classes $[D]$ span $H^2(\P\cH_n)$ and generate
$H^*(\P\cH_n)$ as an algebra.
\end{proposition}

\paragraph{Proof.} The first assertion clearly follows
from the second one. Let us prove the second assertion.
For $n=2$, $\P\cH_2$ is topologically a
sphere, and the class
$[D] = [\{1\} \sqcup \{2\}]$ is a point;
thus the proposition is true. For $n \geq 3$, the classes
$[D]$ corresponding to partitions with at least two
elements in each part
generate the cohomology algebra of $\ocM_n$,
and those corresponding to partitions with a single element
in one part produce the class $\Psi_n$. Therefore, they
generate the cohomology algebra of $\P\cH_n$. \qed

\subsection{Cohomological identities on the Hurwitz spaces}

In this section we study various relations between the
cohomology classes Poincar\'e dual to strata in $\P\cH_n$.
Note that the symmetric group $S_n$ acts on $\cH_n$ by
permuting the poles. Therefore, it also acts on the
cohomology algebra of $\P\cH_n$.
Strata of the stratification, as well as most of other
cohomology classes that
we consider, are invariant under this action.
Therefore,
their cohomology classes belong to the subalgebra
of $S_n$-invariant cohomology classes.

Recall that $C_n$ (respectively $M_n$) denotes the $2$-cohomology class
Poincar\'e dual to the caustic (respectively
the Maxwell stratum) in $\P\cH_n$ (see Example~\ref{ExCM}).

Let $p$ and $q$ be two positive integers,
$p,q \geq 1$, $p+q = n$. We have
introduced the cohomology class $[A \sqcup B]$ assigned to a
partition $A \sqcup B$ of the set 
$\{ 1, \dots, n \}$ into two nonempty parts.
Denote by $\Delta_{p,q}$ the sum of the
classes $[A \sqcup B]$ for all partitions such that
$|A| = p$, $|B| = q$. (If $p=q$ the partitions $A \sqcup B$
and $B \sqcup A$ are two different terms of the sum,
although they determine the same cohomology class.)

\begin{definition}
Denote by $\Delta_n$ the $2$-cohomology class
$$
\Delta_n = \frac12 \,\, \sum_{p+q=n} \Delta_{p,q}
$$
and call it the {\em boundary stratum}.
\end{definition}

The classes introduced above satisfy the following linear relations.

\begin{proposition}\label{Propidentities}
We have
\begin{eqnarray*}
\Psi_n &=& \frac1{2n(n-1)} \,\, \sum_{p+q=n} pq \, \Delta_{p,q}; \\
\Psi_n &=& \frac1{(2n-2)(2n-3)} (3C_n+2M_n+\Delta_n).
\end{eqnarray*}
\end{proposition}

\begin{proposition}\label{Propidentities2}
\begin{eqnarray*}
C_n &=& 6(n-1) \Psi_n - 3 \Delta_n. \\
M_n &=& 2(n-1)(n-6) \Psi_n + 4 \Delta_n. \\
C_n &=& 3 \sum_{p+q=n}
\left( \frac{1}{n}\, pq -\frac12 \right) \Delta_{p,q} \,\, , \\
M_n &=& \sum_{p+q=n}
\left( \frac{n-6}{n}\, pq + 2 \right) \Delta_{p,q} \,\, .
\end{eqnarray*}
\end{proposition}

\paragraph{Proof of Proposition~\ref{Propidentities}.}

The first identity
follows from the identity of Proposition~\ref{PropiDjk},
$$
\Psi_n = [i \D jk],
$$
by summing it over all triples $(i,j,k)$ and regrouping the terms.

In order to prove the second identity
we consider once again the covering $\tPcH_n$ of $\P\cH_n$
introduced in the proof of Theorem~\ref{Thmmult}.
Denote by $\tPsi_n\in H^2(\tPcH_n)$ the pull-back of the class $\Psi_n$
on $\P\cH_n$. In the proof of Theorem~\ref{Thmmult}
we established that
$$
c_1\Bigl(\tLL^*(\cO(-1))\Bigr) = \tPsi_n.
$$
(The pull-back under the lifted Lyashko-Looijenga map of
the class of a hyperplane in $\CP^{2n-3}$ equals $\tPsi_n$.)
Consider a particular union of hyperplanes
in the image space $\CP^{2n-3}$ of $\tLL$.
Namely, the union of all hyperplanes where
the $i$th and the $j$th critical values of $f$
(i.e., the $i$th and the $j$th coordinates
in $\CP^{2n-3}$) are equal (for $1 \leq i < j \leq 2n-2$).
This union represents $(2n-2)(2n-3)/2$ times the class
of a hyperplane. Now consider the preimage of this
union under the map $\tLL$. It consists of the union
of the strata $\widetilde C_n$, $\widetilde M_n$, and
$\widetilde \Delta_n$ (the liftings to $\tPcH_n$ of
the caustic $C_n$, the Maxwell stratum $M_n$, and the
boundary stratum $\Delta_n$). Moreover, $\widetilde M_n$,
$\widetilde \Delta_n$
are simple preimages, but $\widetilde C_n$ is a triple
preimage. (This means that if we take a generic point
$f\in\widetilde C_n$ and its image in $\CP^{2n-3}$,
and then take a generic point close to the image,
then this point will have three preimages close
to $f$.)
Now, the projection of $\tPcH_n$ onto
$\P\cH_n$ is ramified over $C_n$ and $\Delta_n$, and the
order of ramification is~$2$ in both cases.
On the contrary, it is not ramified over $M_n$. This implies
that the geometrical liftings of the strata $C_n$ and
$\Delta_n$ to $\tPcH_n$ represent only one half of the
pull-backs of the corresponding homology classes.
On the other hand, the homology class of the
geometrical lifting of $M_n$ is equal to the lifting
of the homology class of $M_n$. Putting everything
together we obtain
$$
\frac{(2n-2)(2n-3)}2 \,\, \tPsi_n =
3 {\widetilde C_n} + {\widetilde M_n} + {\widetilde \Delta_n},
$$
whence
$$
\Psi_n = \frac{2}{(2n-2)(2n-3)} \,\,
\Bigl(
\frac12 (3C_n) + M_n + \frac12 \Delta_n
\Bigr)
=
\frac1{(2n-2)(2n-3)} \,\, (3C_n+2M_n+\Delta_n).
$$

\qed

\paragraph{Proof of Proposition~\ref{Propidentities2}.}
The proof of this proposition is surprisingly difficult
and uses some results from the subsequent paper~\cite{Zvonkine}.

First of all, note that it suffices to prove the first of
the four identities $C_n = 6(n-1) \Psi_n - 3 \Delta_n.$
The other three identities follow from this one and from
Proposition~\ref{Propidentities}. Moreover, the
$2$-cohomology classes in the left- and the right-hand sides
of the identity are {\em symmetric}, i.e.,
invariant under the action of the
symmetric group $S_n$ by renumbering the poles.
The symmetric part of $H^2(\P\cH_n)$ is spanned
by the classes $\Delta_{p,q}$ (see Section~\ref{Ssec:H*(Hn)}).
Therefore, we must construct a set of $2$-homology
classes of $\P\cH_n$ such that their couplings with the
classes $\Delta_{p,q}$ allow one to distinguish
all linear combinations of the classes $\Delta_{p,q}$.
It is then enough to show that their couplings
with $C_n$ and $6(n-1) \Psi_n - 3 \Delta_n$
coincide. To construct such
a set of $2$-homology classes, we consider several
embeddings of $\CP^1$ into $\P\cH_n$.

For any positive integers $p,q$ such that
$p+q=n$ we are going to consider the following embedding of
$\CP^1$ in $\P\cH_n$. Let $[u:v]$ be the homogeneous
coordinates on $\CP^1$. Consider the following family of
rational functions in the variable $z$:
$$
f_{u,v}(z) = \frac1{uz-a_1} + \dots + \frac1{uz-a_p} +
\frac1{vz-b_1} + \dots + \frac1{vz-b_q},
$$
where the $a_i$ and the $b_j$ are generically chosen complex
numbers. The rational functions of this family do not
belong to $\cH_n$ for $u=0$, $v=0$ or $a_i/u = b_j/v$
(since they have less than $n$ poles); however, there
is a unique way to extend the family to a well-defined
map from $\CP^1$ to $\P\cH_n$, because $\P\cH_n$ is compact.
We denote by $\sigma_{p,q}$ the $2$-homology class of the
image of this map. We are going to study its couplings
with the $2$-cohomology classes 
$C_n$, $\Delta_n$, and $\Delta_{p,q}$.
Note that these classes
are defined as Poincar\'e dual classes
to particular codimension $1$ subvarieties
of $\P\cH_n$. Therefore we will usually speak of their
{\em  intersection indices} (rather than couplings)
with the families $\Delta_{p,q}$.

The three following lemmas describe completely the
intersection indices of the family $\sigma_{p,q}$
with the classes $\Delta_{p',q'}$ and $C_n$.
For shortness, we use the following convention.
If in a lemma we give the intersection indices
$\sigma_{p,q} \cap \Delta_{p',q'}$ and
$\sigma_{p,q} \cap \Delta_{p'',q''}$, then
in the particular case $p'= p''$, $q'= q''$,
the two intersection indices must be added. For example,
when we write
$\sigma_{p,q} \cap \Delta_{p,q} =
\sigma_{p,q} \cap \Delta_{q,p} = 1$,
in the particular case $p=q$ the corresponding
intersection index is equal to $2$.

\begin{lemma} \label{Lem:1}
At each of the points $u=0$ and $v=0$ we have
the following intersection indices:
$$
\sigma_{p,q} \cap \Delta_{p,q} =
\sigma_{p,q} \cap \Delta_{q,p} = 1;
$$
$$
\sigma_{p,q} \cap \Delta_{p',q'}=0
$$
for other pairs $(p',q')$;
$$
\sigma_{p,q} \cap C_n =0.
$$
\end{lemma}

\begin{lemma} \label{Lem:2}
At each of the points $u/v = a_i/b_j$ we have
the following intersection indices:
$$
\sigma_{p,q} \cap \Delta_{1,n-1} =
\sigma_{p,q} \cap \Delta_{n-1,1} = n-2;
$$
$$
\sigma_{p,q} \cap \Delta_{2,n-2} =
\sigma_{p,q} \cap \Delta_{n-2,2} = 1;
$$
$$
\sigma_{p,q} \cap \Delta_{p',q'}=0
$$
for other pairs $(p',q')$;
$$
\sigma_{p,q} \cap C_n =3(n-2).
$$
\end{lemma}

\begin{lemma} \label{Lem:3}
There are $6(pq-1)$ more points of simple intersection
of the family $\sigma_{p,q}$ with $C_n$ corresponding
to rational functions with a double critical point,
defined on a one-component curve.
\end{lemma}

Before proving the lemmas, let us make use of them
to compute the total intersection
indices. (For $\Psi_n$ we use the expression from
Proposition~\ref{Propidentities}.) We have
$$
\sigma_{p,q} \cap \Delta_{1,n-1}=
\sigma_{p,q} \cap \Delta_{n-1,1}= pq(n-2);
$$
$$
\sigma_{p,q} \cap \Delta_{2,n-2}=
\sigma_{p,q} \cap \Delta_{n-2,2}= pq;
$$
$$
\sigma_{p,q} \cap \Delta_{p,q}=
\sigma_{p,q} \cap \Delta_{q,p}= 2;
$$
$$
\sigma_{p,q} \cap \Delta_{p',q'}=0
$$
for other pairs $(p',q')$. And
$$
\sigma_{p,q} \cap \Delta_n = pq(n-1)+2, \qquad
\sigma_{p,q} \cap \Psi_n = pq, \qquad
\sigma_{p,q} \cap C_n = 3pq(n-1)-6.
$$
This is consistent with the identity
$C_n = 6(n-1) \Psi_n - 3 \Delta_n$ that we must prove.

Moreover, looking at the intersection indices with the
classes $\Delta_{p,q}$ we conclude immediately that
the families $\sigma_{p,q}$ for different $p$ and $q$
allow one to distinguish any linear combinations
of the classes $\Delta_{p,q}$. Therefore, the identity
follows from the three lemmas, that we will now prove.

\begin{figure}
\begin{center}
\
\epsfbox{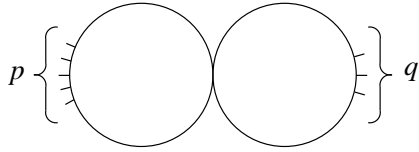}

\caption{\label{Fig:case1} When $u=0$ or $v=0$ the
function is defined on a $2$-component curve with $p$
poles on one component and $q$ on the other one.}
\end{center}
\end{figure}

\paragraph{Proof of Lemma~\ref{Lem:1}.}
Consider the case $u=0$ (the case $v=0$ is similar).
When $u=0$ the function $f$ is defined on a
$2$-component curve shown in Figure~\ref{Fig:case1}.

Such a stable function belongs to the stratum
$\Delta_{p,q} = \Delta_{q,p}$, but not to other
strata $\Delta_{p',q'}$, nor to the caustic $C_n$.
It is easy to check that the intersection with
$\Delta_{p,q} = \Delta_{q,p}$ is simple (unless
$p=q$, in which case it is double by the definition
of $\Delta_{p,p}$). \qed

\begin{figure}
\begin{center}
\
\epsfbox{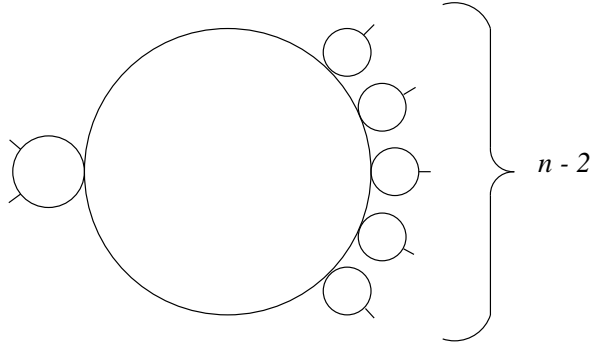}

\caption{\label{Fig:case2} When $u/v = a_i/b_j$ we must
multiply the family by $ub_j - va_i$ in order to get
a finite limit in $\cH_n$ (the limit in $\P\cH_n$ does
not change). Thus the principle part of $n-2$ poles
tend simultaneously to zero, which explains that there
appear $n-2$ small spheres attached to the central 
component.}
\end{center}
\end{figure}

\paragraph{Proof of Lemma~\ref{Lem:2}.}
When $u/v = a_i/b_j$ the function $f$ is
defined on a curve with $n$ irreducible components shown
in Figure~\ref{Fig:case2}. Such a stable function
belongs to $n-2$ different irreducible components
of $\Delta_{1,n-1} = \Delta_{n-1,1}$ and to one
irreducible component of $\Delta_{2,n-2} = \Delta_{n-2,2}$.
It is easy to check that the intersection with each component
is simple.

Now we must find the index of intersection of the
family $\sigma_{p,q}$ with the caustic $C_n$ at the point that
we are considering. This is more difficult, because
this point does not belong to the smooth part of the
caustic. In order to find the intersection index we have to
use a result of the subsequent paper~\cite{Zvonkine},
Theorem~3.9. In the particular case under consideration
the assertion of this theorem is as follows. Consider the
subvariety of $\P\cH_n$ obtained from the stable function
$f$ by moving in all possible ways the $n-1$ points
of intersection of the peripheral components with the
central component of the curve in Figure~\ref{Fig:case2},
without changing the restrictions of $f$ to the
peripheral components. The closure of this subvariety
is isomorphic to $\ocM_{n-1}$. Now, the theorem says that
the neighborhood of this subvariety in $\P\cH_n$ is
isomorphic to the neighborhood of $\ocM_{n-1} \times \{ 0 \}$
in $\cH_{n-1} \times \C^2$. Moreover, the intersection
of the caustic $C_n$ with this neighborhood is isomorphic
to $C_{n-1} \times \C^2$, where $C_{n-1}$ is the
(non-projectivized) caustic in the smaller Hurwitz space.
The intersection of the family $\sigma_{p,q}$ with the
neighborhood is a generic complex curve that intersects the subvariety
$\ocM_{n-1}$ transversally and is not contained
in $C_{n-1} \times \C^2$. Now we claim that the intersection
of the caustic $C_{n-1}$ with a generic fiber of $\cH_{n-1}$
over $\ocM_{n-1}$ is a hypersurface of degree $3(n-3)$, which
is therefore the index of intersection between $\sigma_{p,q}$
and $C_n$ at the point under consideration.

It remains to check that $3(n-3)$ is indeed the degree of the
intersection of $C_{n-1}$ with a generic fiber of $\cH_{n-1}$
over $\ocM_{n-1}$. To do that it suffices to find the degree
of
$$
\frac1{c_1 \dots c_{n-1}}
\mbox{ discrim}_z \, \left( \mbox{numer } \left(
\frac{d}{dz}
\left(
\frac{c_1}{z-z_1} + \dots + \frac{c_{n-1}}{z-z_{n-1}}
\right)
\right)
\right)
$$
as a polynomial in variables $c_1, \dots, c_{n-1}$
(discrim means discriminant, numer means numerator). Indeed,
this polynomial is precisely the equation of the caustic
$C_{n-1}$ in the fiber, where the poles are fixed
at the points $z_1, \dots, z_{n-1}$. The division
by $c_1 \dots c_{n-1}$ is needed because, as one can
easily check, the hyperplanes $c_i=0$ are simple zeroes
of the discriminant, but do not belong to the caustic.
The degree of the above polynomial is
indeed equal to $3(n-3)$ for homogeneity reasons.
\qed

\paragraph{Proof of Lemma~\ref{Lem:3}.}
In order to calculate the intersection of $C_n$ with the
family $\sigma_{p,q}$ outside the points where some
poles get glued together, we will consider the discriminant
of the derivative of the function $f$ of our family
$$
\mbox{ discrim}_z \, \left( \mbox{numer } \left(
\frac{d}{dz}
\left(
\frac1{uz-a_1} + \dots +\frac1{uz-a_p}+
\frac1{vz-b_1} + \dots +\frac1{vz-b_q}
\right)
\right)
\right)
\, .
$$
The degree of this discriminant as a homogeneous polynomial
in $u$ and $v$ equals $2n(2n-3)$. One can check that this
polynomial has double zeroes at the points $u/v = a_i/b_j$.
Moreover, it has a zero of multiplicity $4p^2-6p+3$ at $u=0$
and a zero of multiplicity $4q^2-6q+3$ at $v=0$. Subtracting
the multiplicities of all these zeroes, we obtain
$6(pq-1)$ zeroes outside
the points $u=0, v=0, u/v = a_i/b_j$. \qed

Thus the three lemmas are proved, which completes the proof
of Proposition~\ref{Propidentities2}.

\end{document}